\documentstyle[11pt,epsf]{article}
\newlength{\mytopmargin}
\newlength{\myleftmargin}
\def\zz{\relax\hbox{\small \sf Z\kern-.4em Z}}

\newcommand{\al}{\alpha}
\newtheorem{lemma}{Lemma}[section]
\newtheorem{prop}[lemma]{Proposition}
\begin{document}

\vspace{1cm}
\noindent
\begin{center}{ \Large \bf
SYMMETRIC JACK POLYNOMIALS FROM NON--SYMMETRIC THEORY}
\end{center}
\vspace{5mm}

\begin{center}
T.H.~Baker\footnote{email: tbaker@maths.mu.oz.au; supported by the ARC}
and P.J.~Forrester\footnote{email: matpjf@maths.mu.oz.au; supported by 
the ARC}\\[2mm]
\it Department of Mathematics, University of Melbourne, \\
Parkville, Victoria 3052, Australia
\end{center}
\vspace{.5cm}

\small
\begin{quote} The theory of non-symmetric Jack polynomials is developed
independently of the theory of symmetric Jack polynomials, and this theory
together with the relationship between the non-symmetric, symmetric
and anti-symmetric Jack polynomials is used to deduce the corresponding 
results for the symmetric Jack polynomials. 
\end{quote}

\vspace{.5cm}
\noindent
\section{Introduction}
\setcounter{equation}{0}
The symmetric Jack polynomial $P_\kappa := P_\kappa^{(\alpha)}(z)$,
which are functions of $N$ variables
$z=(z_1,\dots,z_N)$ and labeled by a partition $\kappa$ of length $\le N$, 
can be defined as the unique symmetric polynomial eigenfunction of the
differential operator
\begin{equation}\label{sjeo}
D_2(\alpha) := \sum_{j=1}^N z_j^2 {\partial^2 \over \partial z_j^2}
+ {2 \over \alpha} \sum_{j,k = 1 \atop j \ne k}^N
{z_j^2 \over z_j - z_k}{\partial \over \partial z_j}
\end{equation}
which is of the form
\begin{equation}\label{jex}
P_\kappa^{(\alpha)}(z) = m_\kappa(z) +
\sum_{\mu < \kappa} u_{\kappa \mu} m_\mu(z)
\end{equation}
In (\ref{jex}) $m_\kappa(z)$ is the monomial symmetric function in the
variables $z_1,\dots,z_N$, and the sum is over all partitions $\mu$ which
have the same modulus as $\kappa$ but are smaller in dominance ordering.
The polynomials $P_\kappa$ possess a host of special properties, and in
fact form the natural basis for a class of symmetric
multivariable orthogonal polynomials
generalizing the classical orthogonal polynomials 
\cite{lass91c,macunp1,forr96a}.

Although through the efforts of Macdonald \cite{mac}, Stanley \cite{stan89a}
and others, the theory of symmetric Jack polynomials is highly developed, many
theorems seem difficult to prove. One reason for this is that the symmetric
Jack polynomials are not the most fundamental polynomials in the theory --
this title belongs to the non-symmetric Jack polynomials
$E_\eta := E_\eta(z;\alpha)$, which were introduced \cite{opdam95a}
after the pioneering
works of Macdonald and Stanley. For a given composition of
non-negative integers $\eta := (\eta_1,\dots,\eta_N)$, the polynomials
$E_\eta$ can be defined as the unique polynomial of the form
\begin{equation}\label{estr}
E_\eta(z;\alpha) = z^\eta + \sum_{\nu < \eta} a_{\eta \nu}z^\nu,
\end{equation}
which is an eigenfunction of each of the Cherednik operators
\begin{equation}\label{defch}
\xi_i := \alpha z_i {\partial \over \partial z_i} + \sum_{p < i}
{z_i \over z_i - z_p}(1 - s_{ip}) +
\sum_{p > i} {z_p \over z_i - z_p}(1 - s_{ip}) + 1 - i, \: \: 
(i=1,\dots,N).
\end{equation}
For future reference, we note that the corresponding eigenvalues
are easily calculated as
\begin{equation}\label{evns}
\bar{\eta}_j = \alpha \eta_j - \#\{k<j|\eta_k \ge \eta_j\} -
\#\{k>j|\eta_k > \eta_j\}.
\end{equation}
In (\ref{estr}), $z^\eta =
z_1^{\eta_1} \cdots z_N^{\eta_N}$ and 
$<$ is the partial order on compositions which have the same
modulus, defined by the statement that $\nu < \eta$ if $\nu^+ < \eta^+$
(the superscript $+$ denotes the corresponding partition) with dominance
ordering, or in the case $\nu^+ = \eta^+$, if
$\sum_{j=1}^p \nu_j \le \sum_j^p \eta_j$ for all $p =1,\dots,N$,
while in (\ref{defch}) $s_{ip}$ denotes the transposition operator
which interchanges coordinates $z_i$ and $z_p$.

In this paper we take the view that the non-symmetric Jack polynomials are
more fundamental than the symmetric Jack polynomials, so it should be
possible to first develop the properties of non-symmetric Jack polynomials
independent of the theory of symmetric Jack polynomials, then to use this
theory to develop the properties of the symmetric Jack polynomials. Parts
of this program are already available in the existing literature.
When this is the case we will merely state the result and
give references. However, there are other instances where existing
symmetric Jack polynomial theory has been used to deduce properties of the
non-symmetric Jack polynomials \cite{sahi96a}, as well as cases where the
pathway from the non-symmetric to symmetric theory has yet to be specified.

\section{Non-symmetric Jack polynomial theory}
\setcounter{equation}{0}
The following orthogonality properties are well known and simple to deduce
directly.

\begin{prop}\label{p1} {\rm \cite{opdam95a}}
Define an inner product by
\begin{equation}\label{innerct}
\langle f , g \rangle_I := \int_{[0,1]^N} 
\prod \nolimits_{1 \le j < k \le N}|z_k - z_j|^{2/\alpha }
f(z_1^*,\dots,z_N^*)
g(z_1,\dots,z_N) \, dx_1 \dots dx_N,
\end{equation}
where $z_j := e^{2 \pi i x_j}$ and $*$ denotes the complex conjugate. The
polynomials $\{E_\eta(z;\alpha)\}$ form an orthogonal set with respect to
this inner product.
\end{prop}

\begin{prop}\label{p2} {\rm \cite{sahi96a}}
Define the polynomials $\{q_\eta(x)\}$ by \cite{dunkl96a}
\begin{equation}\label{defq}
\Omega(x,y) = \sum_{\eta} q_\eta(x) y^\eta, \qquad
\Omega(x,y) := {1 \over \prod_{j=1}^N (1 - x_j y_j)}
{1 \over \prod_{j,k=1}^N (1 - x_j y_k)^{1/\alpha}} 
\end{equation}
and an inner product by $\langle q_{\nu}(x) , x^{\eta} \rangle_q := 
\delta_{\nu,\eta}$. Then
\begin{equation}\label{u}
\Omega(x,y) = \sum_{\eta}{1 \over u_\eta} E_\eta(x;\alpha) E_\eta(y;\alpha),
\qquad u_\eta := \langle E_\eta, E_\eta \rangle_q,
\end{equation}
or equivalently the 
polynomials $\{E_\eta(z;\alpha)\}$ form an orthogonal set with respect to
the inner product $\langle \cdot , \cdot \rangle_q$.
\end{prop}

Now, the non-symmetric Jack polynomials can be computed recursively by
using just two types of operators. The first, introduced by Knop and Sahi
\cite{knop96c}, is the raising-type operator $\Phi$, defined to act on 
functions by
\begin{equation}
\Phi f(z_1,\dots,z_N) = z_N f(z_N,z_1,\dots,z_{N-1}),
\end{equation}
and on compositions by
\begin{equation}
\Phi \eta = (\eta_2,\dots,\eta_{N},\eta_1+1),
\end{equation}
while the second is the elementary transposition $s_i := s_{i \, i+1}$
which acts on functions by interchanging the coordinates $z_i$ and 
$z_{i+1}$, and acts on compositions by interchanging the $i$th and $(i+1)$th
parts. On the non-symmetric Jack polynomials these operators have the action
\begin{equation}
\Phi E_\eta = E_{\Phi \eta}
\end{equation}
and
\begin{equation}
s_i E_\eta = \left \{ \begin{array}{ll}
{1 \over \delta_{i,\eta}} E_\eta + (1 - {1 \over \delta^2_{i,\eta}})
E_{s_i\eta} & \eta_i > \eta_{i+1}\\[.1cm]
E_\eta & \eta_i = \eta_{i+1}\\[.1cm]
{1 \over \delta_{i,\eta}} E_\eta + E_{s_i \eta} & \eta_i < \eta_{i+1}
\end{array} \right.
\end{equation}
where $\delta_{i,\eta} := \bar{\eta}_i - \bar{\eta}_{i+1}$.

Using the operators $\Phi$ and $s_i$, 
it is very simple to establish by recurrence formulas
for  $E_\eta(1^N)$
\cite{sahi96a} and $\langle E_\eta, E_\eta \rangle_I$ \cite{forr96d}. 
To write down these formulas requires some notation.
Following Sahi \cite{sahi96a}, for a node $s=(i,j)$ in the diagram of
a composition
define 
the arm length $a(s)$, arm colength $a'(s)$, leg length
$l(s)$ and  leg colength $l'(s)$ by
\begin{eqnarray}
a(s)= \eta_i - j && l(s) = \#\{k>i|j\leq \eta_k\leq\eta_i\} \;+\;
\#\{k<i|j\leq \eta_k+1\leq\eta_i\} \nonumber\\
a'(s)=j - 1 && l'(s) = \#\{k>i| \eta_k > \eta_i\} \;+\;
\#\{k<i|\eta_k\geq\eta_i\}  \label{guion}
\end{eqnarray}
Using these, define constants 
\begin{eqnarray}
d_{\eta} &:=& \prod_{s\in\eta} (\al(a(s)+1) + l(s)+1)\hspace{2cm} 
d'_{\eta} := \prod_{s\in\eta} (\al(a(s)+1) + l(s)) \nonumber\\
e_{\eta} &:=& \prod_{s\in\eta} (\al(a'(s)+1) + N -l'(s)) \hspace{1.7cm} 
e_{\eta}' := \prod_{s\in\eta} (\al(a'(s)+1) +  N -1-l'(s)) \nonumber \\
b_\eta & := & \prod_{s\in\eta} (\al a'(s) + N -l'(s)).
\label{cajamarca}
\end{eqnarray}
We remark that with the generalized factorial defined by
$$
[u]_\kappa^{(\alpha)} := \prod_{j=1}^N 
{\Gamma(u-{1 \over \alpha}(j-1) + \kappa_j) \over \Gamma(1
-{1 \over \alpha}(j-1))},
$$
we have
\begin{equation}\label{e}
e_{\eta} = \alpha^{|\eta|} [1+N/\alpha]^{(\alpha)}_{\eta^+}, \quad
e_{\eta}' = \alpha^{|\eta|} [1+(N-1)/\alpha]^{(\alpha)}_{\eta^+}, \quad
b_\eta = \alpha^{|\eta|} [N/\alpha]^{(\alpha)}_{\eta^+}.
\end{equation}

\begin{prop}\label{p3}{\rm \cite{sahi96a}}
Denoting the non-symmetric Jack polynomial with all variables
$z_j$ set equal to 1 by $E_\eta(1^N)$, we have
\begin{equation}
E_\eta(1^N) = {e_\eta \over d_\eta}.
\end{equation}
\end{prop}

\begin{prop}\label{p4}
Write $\langle E_\eta, E_\eta \rangle_I =: {\cal N}_\eta^{(E)}$. We have
\begin{equation}\label{norml}
{{\cal N}_\eta^{(E)} \over {\cal N}_0^{(E)}} =
{d_\eta' e_\eta \over d_\eta e_\eta'}.
\end{equation}
\end{prop}

\noindent
{\it Remark. } In ref.~\cite[Prop.~2.4]{forr96d} ${\cal N}_\eta^{(E)}$ was
evaluated by recurrence for $1/\alpha = k$, $k \in \zz^+$ only, and
in a form different to (\ref{norml}); however the result (\ref{norml})
was also presented in \cite[eq.~(2.21)]{forr96d} using a method which assumes
the generalized binomial theorem from the theory of symmetric Jack
polynomials \cite{stan89a}. In fact (\ref{norml}) can be established
straightforwardly by verifying the recurrences presented in
\cite[proof of Prop.~2.4]{forr96d}.

\vspace{.2cm}
Next we will show that the generalized binomial theorem can be
deduced within the non-symmetric framework by making use of 
Propositions \ref{p2} and \ref{p4}. This, combined with the
evaluation of ${\cal N}_\eta^{(E)}$ given by (\ref{norml}), allows
the value of $u_\eta$ in (\ref{u}) to be computed.
To derive the generalized binomial theorem, first one replaces
$N$ by $kN$, $k \in \zz^+$ in (\ref{u}), then sets
$y_1,\dots,y_{kN}$ equal to 1 and $x_{N+1},\dots,x_{kN}$ equal to 0.
Noting that $u_\eta$ and $d_\eta$ are independent of $N$, and making
use of the formula in (\ref{e}) for $e_\eta$, gives the formula
\begin{eqnarray}\label{bi1}
\prod_{j=1}^N {1 \over (1-x_j)^{kN/\alpha +1}} & = & 
\sum_{\eta} { e_\eta  |_{N \mapsto kN} \over u_\eta d_\eta}
E_\eta(x;\alpha) \nonumber \\
& = & \sum_{\eta} { 
\alpha^{|\eta|} [1+kN/\alpha]^{(\alpha)}_{\eta^+} 
\over u_\eta d_\eta}
E_\eta(x;\alpha).
\end{eqnarray}
By inspection each term on the r.h.s.~of (\ref{bi1}) is a polynomial in 
$kN/\alpha$. Also, expanding the l.h.s.~as a power series we see that
each term is a polynomial in $kN/\alpha$. Since both sides are equal for
each $k \in \zz^+$, they must in fact be equal for all (complex)
values of $kN/\alpha =: r-1$. Thus (\ref{bi1}) can be rewritten as
\begin{equation}\label{bi2}
\prod_{j=1}^N {1 \over (1-x_j)^r} =
\sum_{\eta} {  \alpha^{|\eta|} [r]^{(\alpha)}_{\eta^+} 
\over  u_\eta d_\eta}
E_\eta(x;\alpha).
\end{equation}
This is eq.~(2.13) of ref.~\cite{forr96d}
with 
\begin{equation}\label{rep}
d_\eta/f_\eta \mapsto u_\eta/d_\eta.
\end{equation}
Now independent of the particular value of $u_\eta$ the reasoning
of ref.~\cite{forr96d}, which relies on (simple) properties of the $E_\eta$,
remains valid and gives
the integration formula (2.19) of ref.~\cite{forr96d} with the replacement
(\ref{rep}), the asymptotics of which implies that
\begin{equation}
\alpha^{|\eta|}{u_\eta \over d_\eta'} {{\cal N}_\eta^{(E)} \over
{\cal N}_0^{(E)}} = {e_\eta \over d_\eta e_\eta'}
\end{equation}
(eq.~(2.21) of ref.~\cite{forr96d} with the replacement (\ref{rep})).
Comparison with (\ref{norml}) allows $u_\eta$ to be specified.

\begin{prop}\label{p5}
We have
\begin{equation}\label{ue}
u_\eta := \langle E_\eta, E_\eta \rangle_q = {d_\eta' \over d_\eta}.
\end{equation}
\end{prop}
{\it Remark. } In ref \cite{sahi96a} (\ref{ue}) was derived by making use 
of results from the theory of symmetric Jack polynomials.

\section{Symmetric Jack polynomial theory}
\setcounter{equation}{0}
In this section we will provide the analogues of Propositions 2.1--2.5 for
the symmetric Jack polynomials, by using the relationship between the
symmetric, anti-symmetric
and non-symmetric Jack polynomials. First, it is well known 
(see e.g.~\cite{forr96c}) how to 
use the Cherednik operators $\xi_i$ to form an operator which separates
the eigenvalues of the $P^{(\alpha)}_{\kappa}(x)$ and is self adjoint 
with respect to
(\ref{innerct}) thus establishing the analogue of Proposition \ref{p1}.

\begin{prop}
The symmetric Jack polynomials $\{P_\kappa^{(\alpha)}(x)\}$ form an
orthogonal set with respect to the inner product (\ref{innerct}).
\end{prop}

To our knowledge there is no previous literature on deducing the
symmetric analogue of Proposition \ref{p2} using the relationship
between the symmetric and non-symmetric Jack polynomials. We will do
this by  applying the antisymmetrization operation Asym, where
$$
{\rm Asym}\, f(x_1,\dots,x_N) :=
\sum_{\sigma \in S_N} (-1)^{\ell(\sigma)} f(x_{\sigma(1)}, \dots,
x_{\sigma(N)}),
$$
to both sides of (\ref{u}). Also required is the formula \cite{forr96c}
\begin{equation}\label{asyme}
{\rm Asym}\, E_\rho(x;\alpha)  = c_\rho \Delta (x) P_{\eta^+}^{
(\alpha/(\alpha+1))}(x), \qquad \Delta (x) := \prod_{1 \le j < k \le N}
(x_j - x_k),
\end{equation}
where $\delta =: (N-1,N-2,\dots,0)$ and all parts of $\rho:=\eta+
\delta$ are assumed distinct
Asym$E_{\rho}(x;\alpha)=0$).
\begin{prop}
We have
\begin{equation}\label{pi}
\Pi^{(\alpha/(\alpha + 1))}(x,y) =
\sum_\kappa {1 \over v_\kappa^{(\alpha/(\alpha+1))}}
P_\kappa^{(\alpha/(\alpha+1))}(x)
P_\kappa^{(\alpha/(\alpha+1))}(y), \quad 
\Pi^{(\alpha)}(x,y) := {1 \over \prod_{j,k=1}^N (1 - x_j y_k)^{1/\alpha}}
\end{equation}
for some constants $v_\kappa^{((\alpha/(\alpha+1))}$ independent of $N$.
Hence, defining the polynomials $g_\kappa(x)$ by
$$
\Pi^{(\alpha)}(x,y) = \sum_{\kappa} g_\kappa(x) m_\kappa(y),
$$
and an inner product by $\langle g_\mu(x), m_\kappa(x) \rangle_g = 
\delta_{\mu,\kappa}$, we have that
\begin{equation}\label{v}
\langle P_\kappa^{(\alpha)}, P_\mu^{(\alpha)} \rangle_g = 
v_\kappa^{(\alpha)} \delta_{\mu,\kappa}.
\end{equation}
\end{prop}
{\bf Proof.} \quad The second statement follows from (\ref{pi}) by a
standard argument \cite{mac}. To derive (\ref{pi}), we apply 
Asym$^{(y)}$ to both sides of (\ref{u}). On the l.h.s.~, according to
the Cauchy double alternant formula we have
\begin{equation}
{\rm Asym}^{(x)} {1 \over \prod_{j=1}^N (1 - x_j y_j)} =
{\Delta(x) \Delta (y) \over \prod_{j,k=1}^N (1 - x_j y_k)},
\end{equation}
and thus
\begin{equation}
{\rm Asym}^{(x)} \Omega (x,y) = {\Delta(x) \Delta(y) \over
\prod_{j,k=1}^N (1 - x_j y_k)^{(\alpha + 1)/\alpha}},
\end{equation}
while on the r.h.s.~we use (\ref{asyme}). This gives
\begin{equation}\label{dd}
{\Delta(x) \Delta(y) \over
\prod_{j,k=1}^N (1 - x_j y_k)^{(\alpha + 1)/\alpha}} =
\Delta(x) \sum_{\rho}{}^*  
{c_\rho \over u_\rho}  P_{\eta^+}^{(\alpha/(\alpha+1))}(x)
E_\rho(y;\alpha),
\end{equation}
where the $*$ denotes that the sum is restricted to $\rho$ with distinct
parts. 

Now applying ${\rm Asym}^{(y)}$ to both sides of (\ref{dd}) we see that
the l.h.s.~is simply multiplied by $N!$, while on the r.h.s.~$E_\eta(y;
\alpha)$ is replaced by using (\ref{asyme}). Cancelling $\Delta(x)
\Delta(y)$ from both sides and summing over permutations of $\eta$ which
give the same $\kappa$ (this only contributes a constant factor to each
term) gives (\ref{pi}). The stability property of the symmetric
Jack polynomials, $P^{(\alpha)}_{\kappa}(z_1,\dots,z_{N-1},0) =
P^{(\alpha)}_{\kappa}(z_1,\dots,z_{N})$ for $\ell(\kappa) \leq N-1$ used in
(\ref{pi}) shows that the $v^{(\alpha)}_{\kappa}$ are independent 
of $N$.\hfill$\Box$

\vspace{.2cm}
The symmetric analogue of Proposition \ref{p3} can be derived by applying
the symmetrization formula \cite{dunkl97a}
\begin{equation}\label{syme}
{\rm Sym} \, E_{\eta^R}(z;\alpha) = \prod_{j=1}^{\eta_1^+}
f_j! \, P_{\eta^+}^{(\alpha)}(z),
\end{equation}
where $\eta^R:= (\eta_N^+, \eta_{N-1}^+, \dots, \eta_1^+)$ and
the $f_j$ denotes the frequency of the part $j$ in $\eta^+$.
The derivation of (\ref{syme}) given in \cite{dunkl97a} is entirely 
within the framework of non-symmetric Jack polynomial theory.

\begin{prop}
We have
\begin{equation}\label{p1n}
P_{\eta^+}^{(\alpha)}(1^N)  = 
{N! \over  \prod_{j=1}^{\eta_1^+}
f_j!} {e_{\eta^R} \over d_{\eta^R}} 
 =  
{ b_{\eta^+} \over h_{\eta^+}}
\end{equation}
where
\begin{equation} 
h_\kappa := \prod_{s \in \kappa}
(\alpha a(s) + l(s) + 1).
\end{equation}
\end{prop}

\begin{figure} 
\epsfxsize=15cm 
\centerline{\epsfbox{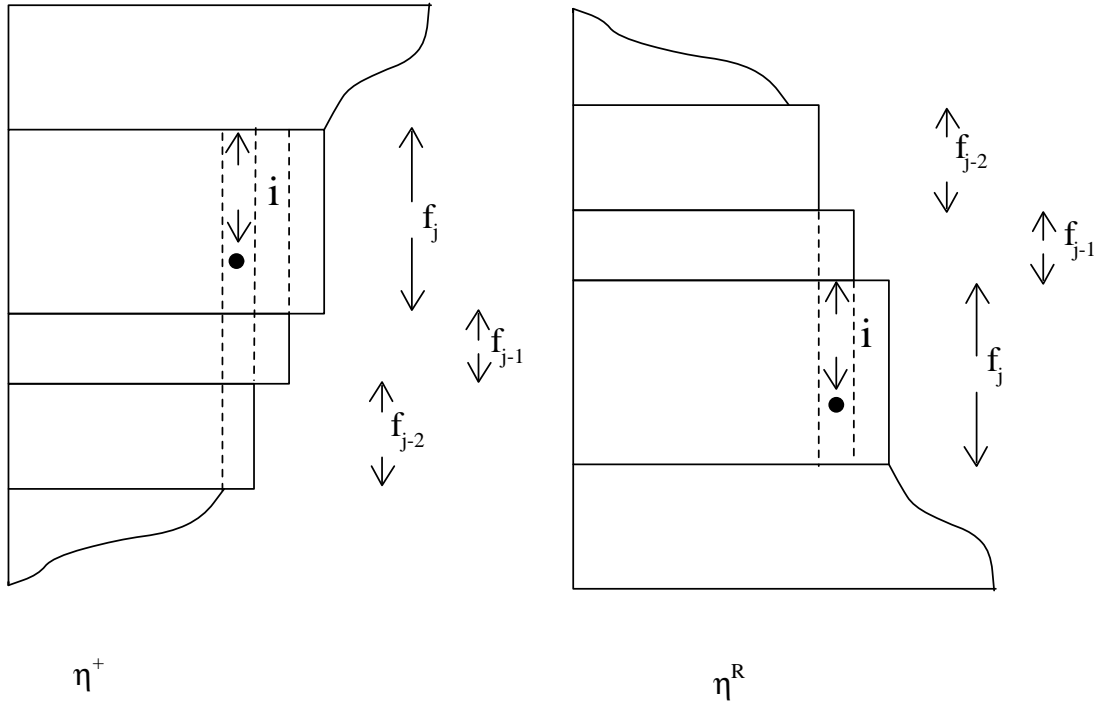}}
\caption{\label{f1} Diagrams specifying the nodes which give the same
contribution to $h_{\eta^+}$ and $d_{\eta^R}$ respectively.}
\end{figure}

\noindent
{\bf Proof. } \quad The first equality is immediate from (\ref{syme}).
For the second equality first note from (\ref{e}) that $e_{\eta^R}
=e_{\eta^+}$ and that
$$
\frac{e_{\eta}}{b_{\eta}} = \frac{e_{\eta^+}}{b_{\eta^+}} =
\frac{1}{N!} \prod_{j=1}^N \left(\alpha \eta_j^+ +N-j+1 \right).
$$
It thus suffices to show that 
\begin{equation} \label{tool}
\frac{h_{\eta^+}}{\prod_i f_i!} = \frac{d_{\eta^R}}
{\prod_{j=1}^N (\alpha \eta_j^+ +N-j+1 )}
\end{equation}
By examining the diagram of $\eta^+$ displaying the repeated parts 
explicitly, the nodes at the {\it end} of each row
contribute $\prod_i f_i!$ to the quantity $h_{\eta^+}$. Similarly, the
nodes at the {\it beginning} of each row of $\eta^R$ contribute 
$\prod_{j=1}^N (\alpha \eta_j^+ +N-j+1 )$ to $d_{\eta^R}$. Moreover, it
can be seen that a node $s\in\eta^+$ (not at the end of a row) in the 
$i$'th row of the $j$'th block has the same leg length as the node
$s'$ in the $i$'th row of the $j$'th block in $\eta^R$ (corresponding
to the node in the $(f_j-i)$'th row in $\eta$), shifted one
column to the right (see fig. \ref{f1}). These nodes 
certainly possess arm lengths
differing by one, and hence give the same contribution to $h_{\eta^+}$
and $d_{\eta^R}$ respectively. Thus (\ref{tool}) is proven. 
\hfill$\Box$

\vspace{.2cm}
Combining (\ref{syme}) with the expansion formula
\begin{equation}\label{pe}
P_\kappa^{(\alpha)}(z) = d_{\eta^+}' \sum_{\eta: \eta^+ = \kappa}
{1 \over d_\eta'} E_\eta(z;\alpha),
\end{equation}
(this was first derived in \cite{sahi96a} using results from the theory of
symmetric Jack polynomials, but, as done in  \cite{dunkl97a}, it
is easily derived within the framework of
non-symmetric Jack polynomial theory), it is simple to
derive the value of $\langle P_\kappa^{(\alpha)}, P_\kappa^{(\alpha)}
\rangle_I =: {\cal N}_\kappa^{(P)} $. 

\begin{prop}
With $\kappa = \eta^+$ we have
\begin{equation} \label{agosto}
{{\cal N}_\kappa^{(P)} \over {\cal N}_0^{(P)}} =
{N! \over \prod_{j=1}^{\eta_1^+}
f_j!} {d_{\eta^+}' e_{\eta^R} \over d_{\eta^R} e_{\eta^R}'} =
{b_{\eta^+} d_{\eta^+}' \over e_{\eta^+}' h_{\eta^+}}.
\end{equation}
\end{prop}

\noindent
{\bf Proof. } \quad The first equality is given in \cite{dunkl97a}, using 
the formulas specified above. The second equality follows from the identity 
implicit in the second equality of (\ref{p1n}).
\hfill$\Box$

\vspace{.2cm}
It remains to establish the analogue of Proposition \ref{p5} and thus
specify the constant $v_\kappa^{(\alpha)}$ in (\ref{v}). This can be done
by proceeding as in the derivation of (\ref{bi2}) and deducing from
(\ref{pi}) and (\ref{p1n}) that
\begin{equation}\label{bi3}
\prod_{j=1}^N {1 \over (1-x_j)^r} =
\sum_{\kappa}
{\alpha^\kappa [r]_\kappa^{(\alpha)} \over v_\kappa^{(\alpha)}  h_\kappa} P_\kappa^{(\alpha)}(x).
\end{equation}
Now substituting (\ref{pe}) for $P_\kappa^{(\alpha)}$ and comparing the 
resulting expression with (\ref{bi2}) with $u_\eta$ therein replaced by
(\ref{ue}) we can read off the value of $v_\kappa^{(\alpha)}$
(a result which  was first
derived in ref.~\cite{stan89a}).

\begin{prop}\label{pv}
We have
\begin{equation}
v_\kappa^{(\alpha)} := 
\langle P_\kappa^{(\alpha)}, P_\kappa^{(\alpha)} \rangle_g =
{ d_\kappa'  \over h_\kappa}.
\end{equation}
\end{prop}

As an extension of the above theory, we will first provide the
evaluation of the constant $c_\rho$ in (\ref{asyme}), and second
investigate the value of the norm of the anti-symmetric Jack polynomial
defined as the monic version of (\ref{asyme}) 
$$
S^{(\alpha)}_{\rho^+}(x) := \Delta(x) P_{\eta^+}^{(\alpha/(\alpha+1))}
(x)
$$
where $\rho^+=\eta^+ + \delta$.

The first task requires the expansion formula
\begin{equation}\label{du}
\Delta (y) P_{\eta^+}^{(\alpha/(\alpha + 1))}(y) =
{(-1)^{N(N-1)/2} \over d_{\rho^+}} 
\sum_{\sigma \in S_N} (-1)^{\ell(\sigma)} d_{\sigma(\rho)}
E_{\sigma(\rho)}(y;\alpha),
\end{equation}
which follows from (\ref{asyme}), and the expansion for
Asym$E_{\eta^R}(y;\alpha)$ given in \cite{dunkl97a}.

\begin{prop}\label{mejorar}
We have
\begin{equation} \label{jeu}
c_\rho = (-1)^{N(N-1)/2} {d'_{\rho}(\alpha) \over d_{\rho^+}(\alpha)}
{h_{\eta^+}(\alpha/(\alpha + 1)) \over 
d_{\eta^+}'(\alpha/(\alpha + 1))}.
\end{equation}
where $\rho^+=\eta^+ + \delta$.
\end{prop}

\noindent
{\bf Proof. } \quad Our strategy is the anti-symmetric analogue of that
used in ref.~\cite{forr96d} to provide the derivation of the constant
$\tilde{c}_\eta$ in the formula Sym$E_\eta(x;\alpha) = \tilde{c}_\eta
P_{\eta^+}^{(\alpha)}(x)$. First substitute the r.h.s.~of
(\ref{pi}) in the l.h.s.~of (\ref{dd}) and cancel the factor of 
$\Delta(x)$ from both sides of the resulting expression. Next
rewrite the sum over $\kappa$ in the new l.h.s.~of (\ref{dd}) as a
sum over distinct partitions $\eta^+$, and substitute
(\ref{du}) for $\Delta(y)  P_{\eta^+}^{(\alpha/(\alpha + 1))}(y)$.
The result follows by equating coefficients of
$P_{\eta^+}^{(\alpha/(\alpha + 1))}(x) E_\rho(y;\alpha)$ on both sides
and using Propositions \ref{p5} and \ref{pv} to substitute for
$u_\rho$ and $v_{\eta^+}^{(\alpha/(\alpha + 1))}$. \hfill $\Box$

It is natural to seek a simplification of the ratio 
$h_{\eta^+}(\alpha/(\alpha + 1)) / d_{\eta^+}'(\alpha/(\alpha + 1))$
in (\ref{jeu}).
In fact, this simplification can be obtained by investigating the value
of the norm of the anti-symmetric Jack polynomial $S_{\rho^+}^{(\alpha)}
(x)$.  It follows from the definition 
of the inner product (\ref{innerct}), and (\ref{agosto}), that
\begin{equation} \label{black}
\langle S_{\rho^+}^{(\alpha)} , S_{\rho^+}^{(\alpha)} \rangle_I^{(\alpha)}
= \langle P_{\eta^+}^{(\alpha/(\alpha+1))} , P_{\eta^+}^{(\alpha/(\alpha+1))} 
\rangle_I^{(\alpha/(\alpha+1))} = 
\frac{b_{\eta^+}(\alpha/(\alpha+1))\,d'_{\eta^+}(\alpha/(\alpha+1))}
{e'_{\eta^+}(\alpha/(\alpha+1))\,h_{\eta^+}(\alpha/(\alpha+1))}\:
{\cal N}_0^{(P)}(\alpha/(\alpha+1)) .
\end{equation}
Here, the superscript in $\langle\cdot , \cdot\rangle_I^{(\zeta)}$
denotes that the weight function is $|\Delta|^{2/\zeta}$. On the other
hand, using the same method as in the derivation of the first equality
in (\ref{agosto}), it was shown in \cite{dunkl97a} that 
\begin{equation} \label{white}
\langle S_{\rho^+}^{(\alpha)} , S_{\rho^+}^{(\alpha)} \rangle_I^{(\alpha)}
= N!\,\frac{d'_{\rho^R}(\alpha)\,e_{\rho^+}(\alpha)}
{d_{\rho^+}(\alpha)\,e'_{\rho^+}(\alpha)}\: {\cal N}_0^{(P)}(\alpha)
\end{equation}
It is our aim in what follows to reconcile these two seemingly different
expressions for the anti-symmetric Jack polynomial norms.

First consider the expression 
$$
(1+\alpha)^{|\eta|} b_{\eta^+}(\alpha/(\alpha+1)) = \prod_{s\in\eta^+}
\left( \alpha(N+a'(s)-l'(s)) + N - l'(s) \right) .
$$
{}From the definitions (\ref{cajamarca}), we see that 
a given node $s\in\eta^+$ contributes the same to this quantity as the
same node $s\in\rho^+/\delta$ (i.e. shifted by $N-i=N-l'(s)-1$ in
row $i$) contributes to $e_{\rho^+}$. Hence
$$
e_{\rho^+}(\alpha) = e_{\delta}(\alpha) (1+\alpha)^{|\eta|}\, 
b_{\eta^+}(\alpha/(\alpha+1)) .
$$
Similarly
$$
e'_{\rho^+}(\alpha) = e'_{\delta}(\alpha) (1+\alpha)^{|\eta|}\, 
e'_{\eta^+}(\alpha/(\alpha+1)) .
$$
Thus
\begin{equation} \label{society.1}
\frac{e_{\rho^+}(\alpha)}{e'_{\rho^+}(\alpha)} = 
\frac{e_{\delta}(\alpha)}{e'_{\delta}(\alpha)}\:
\frac{b_{\eta^+}(\alpha/(\alpha+1))}{e'_{\eta^+}(\alpha/(\alpha+1))}.
\end{equation}

In a similar manner consider 
\begin{equation} \label{skive.1}
(1+\alpha)^{|\eta|} h_{\eta^+}(\alpha/(\alpha+1))
=\prod_{s\in\eta^+}\left( 
\alpha(a(s)+l(s)+1) + l(s)+1\right)
\end{equation}
For every node $s\in\eta^+$ contributing to the quantity (\ref{skive.1})
we need to locate a node $s'\in\rho^+=\eta^+ + \delta$ contributing the
same amount to $d_{\rho^+}(\alpha)$ i.e. we need a node $s'$ such that
\begin{equation} \label{skive.2}
a_{\rho^+}(s') = a_{\eta^+}(s) + l_{\eta^+}(s), \hspace{2cm}
l_{\rho^+}(s') = l_{\eta^+}(s)
\end{equation}

\begin{figure} 
\epsfxsize=15cm 
\centerline{\epsfbox{skive2.ps}}
\caption{\label{f2} Leg-length diagrams of $\eta=(87^243^21)$
and $\rho^+=\eta^+ + \delta$ for $N=9$}
\end{figure}

To see where such nodes lie, it is convenient to visualize $\rho^+$ 
obtained from $\eta^+$ by inserting the {\it columns} of the staircase
partition $\delta$ between the columns of $\eta^+$ such that the column
lengths are weakly decreasing. If one considers the columns of $\delta$
to be placed to the left of the columns of $\eta^+$ of the same length
(these are the shaded columns in fig. \ref{f2}), then it is clear that
the leg lengths of the nodes $s=(i,j)\in\eta^+$ remain unchanged. However
their arm lengths increase by the number of columns of $\delta$ with
column lengths between $i$ and $(\eta^+)_j'-1$ i.e. by
$(\eta^+)_j'-i=l(s)$ ($(\eta^+)'$ denotes the partition conjugate to
$\eta^+$). Thus (\ref{skive.2}) is satisfied and we have
\begin{equation} \label{society.2}
d_{\rho^+}(\alpha) = (1+\alpha)^{|\eta|}\:h_{\eta^+}(\alpha/(\alpha+1))
\; \gamma_{\delta,\eta^+}
\end{equation}
with
$$
\gamma_{\delta,\eta^+} = \prod_{(i,j)\in\delta}\left(\alpha(
a(i,j)+\eta^+_i-\eta^+_{N-j+2}) + l(i,j)+1 \right)
$$

Finally, the diagram of $\rho^R$ is just given by the inverted diagram 
of $\rho^+$. It is clear that the arm length of the nodes remain the same, 
while (due to the definition (\ref{guion})) the leg-length of the inserted
columns of $\delta^R$ increase by one (see fig. \ref{f3}). Similar 
considerations as above imply
\begin{equation} \label{society.3}
d'_{\rho^R}(\alpha) = (1+\alpha)^{|\eta|}\:d'_{\eta^+}(\alpha/(\alpha+1))
\; \gamma_{\delta,\eta^+}
\end{equation}
By using (\ref{society.1}), (\ref{society.2}), (\ref{society.3})
along with the fact that 
$$
\frac{e_{\delta}(\alpha)}{e'_{\delta}(\alpha)}
= \frac{1}{N!} \frac{\prod_{j=1}^N
(j\alpha +N)}{(1+\alpha)^N} = \frac{1}{N!} \frac{{\cal N}^{(P)}_0
(\alpha/(\alpha+1))}{{\cal N}^{(P)}_0(\alpha)}
$$
we complete the reconciliation of (\ref{black}) and (\ref{white}). 

In addition, taking the ratio of (\ref{society.2}) and (\ref{society.3})
provides a simplification of $h_{\eta^+}(\alpha/(\alpha+1)) /
d'_{\eta^+}(\alpha/(\alpha+1))$ and thus a simplification of the 
formula in Proposition \ref{mejorar}. 
\bigskip\\
{\bf Proposition \ref{mejorar}$'$}\quad  
{\it The constant $c_{\rho}$ in
(\ref{asyme}) is given by
$$
c_{\rho} = (-1)^{N(N-1)/2}\: \frac{d'_{\rho}(\alpha)}
{d'_{\rho^R}(\alpha)}
$$
}

\noindent
{\it Remark. } \quad In the case $\rho = \rho^R$, this gives
$c_\rho = (-1)^{N(N-1)/2}$, which is a result contained in
ref.~\cite[lemma 2.1]{dunkl97a}.

\begin{figure} 
\epsfxsize=10cm
\centerline{\epsfbox{skive3.ps}}
\caption{\label{f3} Leg-length diagram for $\rho^R$} 
\end{figure}


\end{document}